\documentclass[12pt]{amsart}
\advance\textheight\topskip
\allowdisplaybreaks

\usepackage{mathabx}
\newcommand{\Smash}{\mathbin{\hash}}
\newcommand{\Braid}{\mathbin{\mbox{\large${\Join}$}}}
\newcommand{\braid}{\mathbin{\mbox{\normalsize${\Join}$}}}

\usepackage{amsfonts,amssymb}
\usepackage{times}
\usepackage[mathscr]{eucal}
\usepackage{mathptmx}
\usepackage{amsthm}

\newcommand{\Dz}{\dd}

\newcommand{\up}[1]{^{{\scriptscriptstyle(#1)}}}

\newcommand{\zero}{_{_{(0)}}}
\newcommand{\mone}{_{_{(-1)}}}

\newcommand{\leftreg}{\kern1pt{\rightharpoonup}\kern1pt}
\newcommand{\rightreg}{\kern1pt{\leftharpoonup}\kern1pt}

\newcommand{\leftregchi}{\kern1pt{\rightharpoonup}_{\!\!\!\!_{\chi}}\kern1pt}

\newcommand{\doteta}{\kern2pt{\cdot}_{\!_{\eta}}\kern1pt}

\newcommand{\Czd}{\oC_{\q}[z,\Dz]}
\newcommand{\Mat}{\mathrm{Mat}}
\renewcommand{\kappa}{\varkappa}
\newcommand{\qfac}[1]{[#1]!\,}
\newcommand{\qint}[1]{[#1]}

\newcommand{\leftii}{\mathbin{\mbox{\footnotesize${\vartriangleright}$}}}

\newcommand{\Dp}[1]{\,_{\phantom{h}}^{\underline{#1\kern-.5pt}\kern.5pt}}

\newcommand{\Sinv}{{S^*}^{-1}}

\newcommand{\Ri}{R\up{1}}
\newcommand{\Rii}{R\up{2}}

\usepackage[bookmarks=false, hyperfigures=false, a4paper]{hyperref}

\newcommand{\id}{\mathrm{id}}
\newcommand{\eval}[2]{\langle#1,\,#2\rangle\,}

\newcommand{\DD}{\mathscr{D}}
\newcommand{\HD}{\mathscr{H}}
\newcommand{\HH}{\boldsymbol{\mathscr{H}}}
\newcommand{\HHH}{\pmb{\textsf{H}}}
\newcommand{\bDDB}{\kern4pt\overline{\kern-3pt\mathscr{D}(B)\kern-3pt}\kern4pt}

\newcommand{\bref}[1]{\textbf{\ref{#1}}}

\renewcommand{\geq}{\geqslant}
\renewcommand{\leq}{\leqslant}

\newcommand{\tensor}{\otimes}

\newcommand{\q}{\mathfrak{q}}

\newcommand{\UresSL}[1]{\overline{\mathscr{U}}_{\q} s\ell(#1)}
\newcommand{\HresSL}[1]{\overline{\mathscr{H}}_{\q} s\ell(#1)}

\newcommand{\mfrac}[2]{\raisebox{.8pt}{\mbox{\small$\displaystyle\frac{#1}{#2}$}}}
\newcommand{\ffrac}[2]{\raisebox{.5pt}{\mbox{\footnotesize$\displaystyle\frac{#1}{#2}$}}}
\newcommand{\fffrac}[2]{\raisebox{.9pt}{\mbox{\scriptsize$\displaystyle\frac{#1}{#2}$}}}
\newcommand{\half}{%
  \mathchoice{\ffrac{1}{2}}{\frac{1}{2}}{\frac{1}{2}}{\frac{1}{2}}}

\newcommand{\qbin}[2]{\mathchoice%
  {\qbinm{#1}{#2}}{\qbinmm{#1}{#2}}%
  {\qbinmm{#1}{#2}}{\qbinmm{#1}{#2}}}
\newcommand{\qbinm}[2]{\mbox{\footnotesize$\displaystyle
    \genfrac{[}{]}{0pt}{}{#1}{#2}$}}
\newcommand{\qbinmm}[2]{\genfrac{[}{]}{0pt}{}{#1}{#2}}

\newcommand{\dd}{\partial}

\newcommand{\oC}{\mathbb{C}}

\newcommand{\oZ}{\mathbb{Z}}

\numberwithin{equation}{section}
\makeatletter
\@addtoreset{equation}{section}
\@addtoreset{subsubsection}{section}

\def\@secnumfont{\bfseries}
\def\subsubsection{\@startsection{subsubsection}{3}%
  \z@{.5\linespacing\@plus.7\linespacing}{-.5em}%
  {\normalfont\bfseries}}
\def\paragraph{\@startsection{paragraph}{4}%
  \z@\z@{-\fontdimen2\font}%
  \normalfont\bfseries}
\def\subparagraph{\@startsection{subparagraph}{5}%
  \z@\z@{-\fontdimen2\font}%
  \normalfont\bfseries}

\makeatother

\swapnumbers
\newtheorem{Thm}[subsection]{Theorem}

\newtheorem{lemma}[subsubsection]{Lemma}

\theoremstyle{definition}
\newtheorem{Dfn}[subsection]{Definition}

\begin{document}

\title{Yetter--Drinfeld structures on Heisenberg doubles and
  chains}

\author[Semikhatov]{A.M.~Semikhatov}

\address{Lebedev Physics Institute
  \hfill\mbox{}\linebreak \texttt{ams@sci.lebedev.ru}}

\begin{abstract}
  For a Hopf algebra $B$ with bijective antipode, we show that the
  Heisenberg double $\HD(B^*)$ is a braided commutative
  Yetter--Drin\-feld module algebra over the Drinfeld double $\DD(B)$.
   The braiding structure allows generalizing $\HD(B^*)\cong
  B^{*\mathrm{cop}}\braid B$ to ``Heisenberg $n$-tuples'' and
  ``chains'' $\dots \braid B^{*\mathrm{cop}}\braid B \braid
  B^{*\mathrm{cop}}\braid B\braid\dots$, all of which are
  Yetter--Drin\-feld $\DD(B)$-module algebras.  For $B$ a particular
  Taft Hopf algebra at a $2p$th root of unity, the construction is
  adapted to yield Yetter--Drin\-feld module algebras over the
  $2p^3$-dimensional quantum group~$\UresSL2$.
\end{abstract}

\maketitle
\thispagestyle{empty}

\section{Introduction}
We establish the properties of $\HD(B^*)$\,---\,the Heisenberg double
of a (dual) Hopf algebra\,---\,relating it to two popular structures:
Yetter--Drin\-feld modules (over the Drinfeld double~$\DD(B)$) and
braiding.  In fact, we construct new examples of Yetter--Drinfeld
\textit{module algebras}, some of which are in addition braided
commutative. 

Heisenberg doubles~\cite{[AF],[RSts],[Sts],[Lu-double]} have been the
subject of some attention, notably in relation to Hopf algebroid
constructions~\cite{[Lu-alg],[P],[BM]} (the basic observation being
that $\HD(B^*)$ is a Hopf algebroid over $B^*$~\cite{[Lu-alg]}) and
also from various other
standpoints~\cite{[K],[Ka],[Mi],[IP]}.\footnote{The ``true,''
  underlying motivation (deriving
  from~\cite{[FGST],[FGST2],[S-q],[G],[S-U],[S-H]}) of our interest in
  $\HD(B^*)$ is entirely left out here.}  We show that they are a rich
source of Yetter--Drinfeld $\DD(B)$-module algebras: $\HD(B^*)$ is a
Yetter--Drin\-feld module algebra over the Drinfeld double $\DD(B)$;
it is, moreover, braided commutative.  Reinterpreting the construction
of $\HD(B^*)$ in terms of the braiding in the Yetter--Drin\-feld
category then allows generalizing Heisenberg \textit{doubles} to
``\textit{$n$-tuples},'' or ``Heisenberg chains''\footnote{A slight
  mockery of the statistical-mechanics meaning of a ``Heisenberg
  chain'' may give way to a genuine, and deep, relation in the context
  of the previous footnote.} (cf.~\cite{[NSz]}), which are all
Yetter--Drin\-feld $\DD(B)$-module
algebras.\enlargethispage{\baselineskip}

In Sec.~\ref{sec:HD}, we establish that $\HD(B^*)$ is a
Yetter--Drin\-feld $\DD(B)$-module algebra, and in
Sec.~\ref{sec:braiding} that it is braided ($\DD(B)$-)
commutative~\cite{[CFM]}; there, $B$ denotes a Hopf algebra with
bijective antipode.  In Sec.~\ref{sec:sl2}, where we construct
Yetter--Drinfeld module algebras for a quantum $s\ell(2)$ at an even
root of unity~\cite{[AGL],[Su],[X],[FGST],[FGST2]}, $B$ becomes a
particular Taft Hopf algebra.

For the left and right regular actions of a Hopf algebra $B$ on $B^*$,
we use the respective notation $b\leftreg
\beta=\eval{\beta''}{b}\beta'$ and $\beta\rightreg
b=\eval{\beta'}{b}\beta''$, where $\beta\in B^*$ and $b\in B$ (and
$\eval{~}{~}$ is the evaluation).  The left and right regular actions
of $B^*$ on $B$ are $\beta\leftreg b=\eval{\beta}{b''}b'$ and
$b\rightreg\beta=\eval{\beta}{b'}b''$.  We assume the precedence
$ab\rightreg\beta=(ab)\rightreg\beta$, $\alpha\beta\leftreg
a=(\alpha\beta)\leftreg a$, and so~on.  For a Hopf algebra $H$ and a
left $H$-comodule $U$, we write the coaction $\delta:U\to H\tensor U$
as $\delta (u)=u\mone\!\tensor u\zero$; then
$\eval{\varepsilon}{\!u\mone}u\zero =u$ and $u\mone'\!\tensor
u\mone''\!\tensor u\zero =u\mone\!\tensor u\zero{}\mone\!\tensor
u\zero{}\zero$.

\section{$\HD(B^*)$ as a Yetter--Drin\-feld $\DD(B)$-module
  algebra}\label{sec:HD}
The purpose of this section is to show that $\HD(B^*)$ is a
Yetter--Drin\-feld $\DD(B)$-module algebra.  The key ingredients are
the $\DD(B)$-comodule algebra structure from~\cite{[Lu-double]}, which
we recall in~\bref{H-comult}, and the $\DD(B)$-module algebra
structure from~\cite{[S-H]}, which we recall in~\bref{H-action}.  The
claim then follows by direct computation.

\subsection{The Heisenberg double $\HD(B^*)$}
The Heisenberg double $\HD(B^*)$ is the smash product $B^*\Smash B$
with respect to the left regular action of $B$ on $B^*$, which means
that the composition in $\HD(B^*)$ is given~by
\begin{equation}\label{H-comp}
  (\alpha\Smash a)(\beta\Smash b)=
  \alpha(a'\leftreg\beta)\Smash a'' b,
  \qquad
  \alpha,\beta\in B^*,\quad
  a,b\in B.
\end{equation}

\subsubsection{}\label{H-comult}
We recall from~\cite{[Lu-double]} that $\HD(B^*)$ can also be obtained
by twisting the product on the Drinfeld double $\DD(B)$ (see
Appendix~\ref{app:D-double}) as follows.  Let
\begin{equation*}
  \eta:\DD(B)\tensor\DD(B)\to k
\end{equation*}
be given by
\begin{equation*}
  \eta(\mu\tensor m,\nu\tensor n)
  =\eval{\mu}{1}\eval{\varepsilon}{n}\eval{\nu}{m}
\end{equation*}
and let $\doteta:\DD(B)\tensor\DD(B)\to\DD(B)$ be defined as
\begin{equation*}
  M\doteta N=M' N' \eta(M'',N''),\qquad M,N\in\DD(B).
\end{equation*}
A simple calculation shows that ${}\doteta{}$ coincides with the
product in~\eqref{H-comp}:
\begin{equation*}
  (\mu\tensor m)\doteta(\nu\tensor n)
  =\mu(m'\leftreg\nu)\tensor m'' n,
  \qquad \mu,\nu\in B^*,\quad m,n\in B.
\end{equation*}

{}From this construction of $\HD(B^*)$, it readily
follows~\cite{[Lu-double]} that the coproduct of $\DD(B)$, viewed as a
map
\begin{equation}\label{delta}
  \begin{split}
    \delta:\HD(B^*)&\to\DD(B)\tensor\HD(B^*)\\
    \beta\Smash b&\mapsto(\beta''\tensor b')
    \tensor
    (\beta'\Smash b''),
  \end{split}
\end{equation}
makes $\HD(B^*)$ into \textit{a left $\DD(B)$-comodule
  algebra} (i.e., $\delta$ is an algebra morphism).

\subsubsection{}\label{H-action}
Simultaneously, \textit{$\HD(B^*)$ is a $\DD(B)$-module algebra},
i.e.,
\begin{equation}\label{m-a}
  M\leftii(A C) = (M'\leftii A)(M''\leftii C)
\end{equation}
for all $M\in\DD(B)$ and $A,C\in\HD(B^*)$, under the $\DD(B)$ action
defined in~\cite{[S-H]}:\pagebreak[3]
\begin{multline}\label{the-action}
  (\mu\tensor m)\leftii(\alpha\Smash a)
  =\mu'''(m'\leftreg\alpha)\Sinv(\mu'')
  \Smash\bigl((m'' a S(m'''))\rightreg\Sinv(\mu')\bigr),\\[-3pt]
  \mu\tensor m\in\DD(B),\quad \alpha\Smash a\in\HD(B^*).
\end{multline}

Evidently, the right-hand side here
factors into the actions of $B^{*\mathrm{cop}}$ and $B$:
\begin{equation*}
  (\mu\tensor m)\leftii(\alpha\Smash a)=
  (\mu\tensor 1)\leftii
  \bigl((\varepsilon\tensor m)\leftii(\alpha\Smash a)\bigr),
\end{equation*}
where
\begin{align*}
  (\varepsilon\tensor m)\leftii(\alpha\Smash a)
  &=(m'\leftreg\alpha)
  \Smash(m'' a S(m'''))\\[-4pt]
  \intertext{and}
  \smash[t]{(\mu\tensor 1)\leftii(\alpha\Smash a)}
  &\smash[t]{{}=\mu'''\alpha\Sinv(\mu'')
    \Smash(a\rightreg\Sinv(\mu'))}.
\end{align*}
This allows verifying that~\eqref{the-action} is indeed an action of
$\DD(B)$ independently of the argument in~\cite{[S-H]}: it suffices to
show that the actions of $B^{*\mathrm{cop}}$ and $B$ taken in the
``reverse'' order combine in accordance with the Drinfeld double
multiplication, i.e., to show that
\begin{align}\label{to-show-action}
  (\varepsilon\tensor m)\leftii
  \bigl((\mu\tensor 1)\leftii(\alpha\Smash a)\bigr)
  &=\bigl((\varepsilon\tensor m)(\mu\tensor 1)\bigr)
  \leftii(\alpha\Smash a)\\
  &=\bigl((m'\leftreg\mu\rightreg S^{-1}(m'''))\tensor m''\bigr)
  \leftii(\alpha\Smash a).\notag
\end{align}
We do this in~\bref{app-prove-Drinfeld}.

The $\DD(B)$-module algebra property was shown in~\cite{[S-H]}, but a
somewhat less bulky proof can be given by considering the actions of
$\mu\tensor 1$ and $\varepsilon\tensor m$ separately.  The routine
calculations are in~\bref{app-prove-module-algebra}.

\begin{Thm}\label{thm:YD}
  $\HD(B^*)$ is a \textup{(}left--left\textup{)} Yetter--Drin\-feld
  $\DD(B)$-module algebra.
\end{Thm}
By this we mean a left module algebra and a left comodule algebra with
the Yetter--Drin\-feld compatibility condition
\begin{equation}\label{eq:YD}
  (M'\leftii A)\mone M''\tensor(M'\leftii A)\zero
  =M' A\mone\tensor(M''\leftii A\zero)
\end{equation}
for all $M\in\DD(B)$ and $A\in\HD(B^*)$.  (For Yetter--Drin\-feld
modules, see~\cite{[Y],[LR],[RT],[Mont],[Sch],[CFM]}.)
Condition~\eqref{eq:YD} has to be shown for the $\DD(B)$ action and
coaction in~\eqref{the-action} and~\eqref{delta}.

\subsubsection{Proof of~\ref{thm:YD}}
To simplify the calculation leading to~\eqref{eq:YD}, we again use
that the action of $\mu\tensor m\in\DD(B)$ factors through the actions
of $\mu\tensor 1$ and $\varepsilon\tensor m$.

First, for $M = \varepsilon\tensor m$, we evaluate the left-hand side
of~\eqref{eq:YD} as
\begin{align*}
  \bigl((\varepsilon\tensor m')
    \leftii(\alpha\Smash a)\bigr)\mone (\varepsilon\tensor m'')
    \tensor
    ((\varepsilon\tensor m')\leftii(\alpha\Smash a))\zero\kern-230pt
  \\
    &=
    \bigl((m\up{1}\leftreg\alpha)'' \tensor (m\up{2} a S(m\up{3}))'
    (\varepsilon\tensor m\up{4})\bigr)
    \tensor
    \bigl(
    (m\up{1}\leftreg\alpha)'\Smash(m\up{2} a S(m\up{3}))''
    \bigr)
    \\
    &=
    \bigl((m\up{1}\leftreg\alpha'') \tensor (m\up{2} a S(m\up{3}))'m\up{4}
    \bigr)
    \tensor
    \bigl(
    \alpha'\Smash(m\up{2} a S(m\up{3}))''
    \bigr)
    \\
    &=
    \bigl((m\up{1}\leftreg\alpha'') \tensor
    m\up{2} a' S(m\up{5})m\up{6}
    \bigr)
    \tensor
    \bigl(
    \alpha'\Smash m\up{3} a'' S(m\up{4})
    \bigr)
    \\
    &=
    \bigl((m\up{1}\leftreg\alpha'') \tensor
    m\up{2} a' \bigr)
    \tensor
    \bigl(
    \alpha'\Smash m\up{3} a'' S(m\up{4})
    \bigr)
\end{align*}
but the right-hand side is given by 
\begin{multline*}
  \mbox{}\\*[-1.3\baselineskip]
  \shoveleft{\bigl((\varepsilon\tensor m')(\alpha''\tensor a')\bigr)
    \tensor
    \bigl((\varepsilon\tensor m'')\leftii(\alpha'\Smash a'')\bigr)}
  \\
  \begin{aligned}[t]
    &=
    \bigl((m\up{1}\leftreg\alpha''\rightreg S^{-1}(m\up{3}))
    \tensor m\up{2} a'\bigr)
    \tensor
    \bigl((m\up{4}\leftreg\alpha')\Smash m\up{5} a'' S(m\up{6})\bigr)
    \\
    &=
    \bigl((m\up{1}\leftreg\alpha'')
    \tensor m\up{2} a'\bigr)
    \tensor
    \bigl((m\up{4}S^{-1}(m\up{3})\leftreg\alpha')
    \Smash m\up{5} a'' S(m\up{6})\bigr)
  \end{aligned}
\end{multline*}
(because 
$\alpha'\tensor(\alpha''\rightreg m)=(m\leftreg\alpha')\tensor\alpha''$),
which is the same as the left-hand side.

Second, for $M = \mu\tensor 1$, using the $\DD(B)$-identity
\begin{equation}\label{identity}
  \bigl(\varepsilon\tensor(a\rightreg\Sinv(\mu''))\bigr)(\mu'\tensor 1)
  = \mu''\tensor(\Sinv(\mu')\leftreg a),
\end{equation}
we evaluate the left-hand side of~\eqref{eq:YD} as
\begin{align*}
  \mbox{}\kern-20pt
  \bigl((\mu''\tensor 1)\leftii(\alpha\Smash a)\bigr)\mone
    (\mu'\tensor 1)
    \tensor
    \bigl((\mu''\tensor 1)\leftii(\alpha\Smash a)\bigr)\zero\kern-270pt
  \\
    &=
    \Bigl(\bigl(\mu\up{4}\alpha\Sinv(\mu\up{3})\bigr)''
    \tensor\bigl(a\rightreg\Sinv(\mu\up{2})\bigr)'
    (\mu\up{1}\tensor 1)\Bigr)
    \\[-2pt]
    &\qquad\qquad\qquad\qquad\qquad\qquad\qquad
    \tensor
    \Bigl(\bigl(\mu\up{4}\alpha\Sinv(\mu\up{3})\bigr)'
    \Smash\bigl(a\rightreg\Sinv(\mu\up{2})\bigr)''\Bigr)
    \\
    &=
    \Bigl(\!\bigl(\mu\up{6}\alpha''\Sinv(\mu\up{3})
    \tensor(a'\rightreg\Sinv(\mu\up{2}))\bigr)
    (\mu\up{1}\tensor 1)\!\Bigr)
    \tensor
    \bigl(\mu\up{5}\alpha'\Sinv(\mu\up{4})
    \Smash a''\bigr)
    \\
    &\stackrel{\eqref{identity}}{=}\!
    \bigl(\mu\up{6}\alpha''\Sinv(\mu\up{3})\mu\up{2}
    \tensor(\Sinv(\mu\up{1})\leftreg a')\bigr)
    \tensor
    \bigl(\mu\up{5}\alpha'\Sinv(\mu\up{4})
    \Smash a''\bigr)
    \\
    &=
    \bigl(\mu\up{4}\alpha''
    \tensor(\Sinv(\mu\up{1})\leftreg a')\bigr)
    \tensor
    \bigl(\mu\up{3}\alpha'\Sinv(\mu\up{2})
    \Smash a''\bigr)
    \\
    &=
    \bigl(\mu\up{4}\alpha''
    \tensor a'\bigr)
    \tensor
    \bigl(\mu\up{3}\alpha'\Sinv(\mu\up{2})
    \Smash (a''\rightreg\Sinv(\mu\up{1}))\bigr)
    \\
    &=\bigl((\mu''\tensor 1)(\alpha''\tensor a')\bigr)
    \tensor
    \bigl((\mu'\tensor 1)\leftii(\alpha'\Smash a'')\bigr),
\end{align*}
which is the right-hand side.

\section{$\HD(B^*)$ as a braided commutative algebra}
\label{sec:braiding}
The category of Yetter--Drin\-feld modules is well known to be braided,
with the braiding $c^{\vphantom{1}}_{U,V}:U\tensor V\to V\tensor U$
given by
\begin{equation*}
  c^{\vphantom{1}}_{U,V}:u\tensor v\mapsto (u\mone\leftii v)\tensor u\zero.
\end{equation*}
The inverse is $c^{-1}_{U,V}:v\tensor u\mapsto u\zero\tensor
S^{-1}(u\mone)\leftii v$.

\begin{Dfn}
  A left $H$-module and left $H$-comodule algebra $X$ is said to be
  \textit{braided commutative}~\cite{[BM]} (or
  $H$-commutative~\cite{[CFM],[CGW]}) if
  \begin{equation}\label{braided-comm}
    y x = (y\mone\leftii x) y\zero
  \end{equation}
  for all $x,y\in X$.
\end{Dfn}

\begin{Thm}\label{thm:braided-comm}
  $\HD(B^*)$ is braided commutative with respect to the braiding
  associated with the Yetter--Drinfeld module structure.
\end{Thm}

\vspace*{-.6\baselineskip}

\subsubsection{Remarks}
\begin{enumerate}
\item The braided$/$$H$-commutativity property may be compared with
  ``quantum commutativity''~\cite{[CW]}.  We recall that for a
  \textit{quasitriangular} Hopf algebra $H$, its module algebra $X$ is
  called quantum commutative if
  \begin{equation}\label{R-comm}
    y x = (\Rii\leftii x)(\Ri\leftii y) \equiv
    {}\cdot{}(R_{21}\leftii(x\tensor y)),\qquad
    x,y\in X,
  \end{equation}
  where $R=\Ri\tensor\Rii\in H\tensor H$ is the universal $R$-matrix
  (and the dot denotes the multiplication in~$X$). A minor source of
  confusion is that this useful property (see,
  e.g.,~\cite{[CW],[Lu-alg],[P]}) is sometimes also referred to as
  $H$-commutativity~\cite{[CW]}.  For a Yetter--Drin\-feld module
  algebra $X$ over a quasitriangular $H$, the properties
  in~\eqref{braided-comm} and \eqref{R-comm} are different (for
  example, a ``quantum commutative'' analogue of
  Theorem~\bref{thm:braided-comm} does not hold for $\HD(B^*)$).  We
  therefore consistently speak of~\eqref{braided-comm} as of ``braided
  commutativity'' (this term is also used in~\cite{[Mj]} in related
  contexts, although in more than one).

\item The two properties, Eqs.~\eqref{braided-comm} and
  \eqref{R-comm}, are ``morally'' similar, however.  To see this,
  recall that a Yetter--Drin\-feld $H$-module is the same thing as a
  $\DD(H)$-module, the $\DD(H)$ action on a left--left
  Yetter--Drin\-feld module~$X$ being defined~as
  \begin{equation*}
    \quad(p\tensor h)\leftii x
    = \eval{\Sinv(p)}{(h\leftii x)\mone}(h\leftii x)\zero,
    \qquad p\in H^*,\quad h\in H,\quad x\in X.
  \end{equation*}
  Let then 
  \begin{equation*}
    \mathscr{R} = \sum_A(\varepsilon \tensor e_A)\tensor (e^A\tensor 1)
    \in\DD(H)\tensor\DD(H)
  \end{equation*}
  be the universal $R$-matrix for the double.  It follows that
  \begin{multline*}
    \qquad{}\cdot{}(\mathscr{R}^{-1}\leftii(x\tensor y))
    =\bigl((\varepsilon \tensor S(e_A))\leftii x\bigr)
    \bigl((e^A\tensor 1)\leftii y\bigr)
    \\
    {}=\eval{e^A}{S^{-1}(y\mone)}
    \bigl(S(e_A)\leftii x\bigr)y\zero
    =(y\mone\leftii x)y\zero
  \end{multline*}
  for all $x,y\in X$, and therefore the braided commutativity property
  can be equivalently stated in the form
  \begin{equation*}
    y x =
    {}\cdot{}(\mathscr{R}^{-1}\leftii(x\tensor y))
  \end{equation*}
  similar to Eq.~\eqref{R-comm} (the occurrence of $\mathscr{R}^{-1}$
  instead of $\mathscr{R}_{21}$ may be attributed to our choice of
  left--left Yetter--Drin\-feld modules).
\end{enumerate}

\subsubsection{Proof of~\bref{thm:braided-comm}}
We evaluate the right-hand side of~\eqref{braided-comm} for
$X=\HD(B^*)$ as
\begin{align*}
  \bigl((\beta\Smash b)\mone\leftii(\alpha\Smash a)\bigr)
    (\beta\Smash b)\zero\kern-120pt\\
    &=((\beta''\tensor b')\leftii(\alpha\Smash a))
    (\beta'\Smash b'')\\
    &=\Bigr(
    \beta\up{4}(b\up{1}\leftreg\alpha)\Sinv(\beta\up{3})
    \Smash \bigl(b\up{2} a S(b\up{3})\rightreg\Sinv(\beta\up{2})\bigr)
    \Bigl)
    (\beta\up{1}\Smash b\up{4})\\
    &=
    \bigl(\beta\up{4}(b\up{1}\leftreg\alpha)\Sinv(\beta\up{3})
    \bigl(b\up{2} a S(b\up{3})\rightreg\Sinv(\beta\up{2})\bigr)'
    \leftreg\beta\up{1}\bigr)\\[-2pt]
    &\qquad\qquad\qquad\qquad\qquad\qquad\qquad\qquad\qquad
    \Smash
    \bigl(b\up{2} a S(b\up{3})\rightreg\Sinv(\beta\up{2})\bigr)''b\up{4}
    \\
    &=
    \bigl(\beta\up{4}(b\up{1}\leftreg\alpha)\Sinv(\beta\up{3})
    \bigl((b\up{2} a S(b\up{3}))'\rightreg\Sinv(\beta\up{2})\bigr)
    \leftreg\beta\up{1}\bigr)\\[-2pt]
    &\qquad\qquad\qquad\qquad\qquad\qquad\qquad\qquad\qquad
    \Smash
    (b\up{2} a S(b\up{3}))''b\up{4}
    \\[-2pt]
    &\stackrel{\checkmark}{=}
    \beta\up{5}(b\up{1}\leftreg\alpha)\Sinv(\beta\up{4})
    \beta\up{1}\eval{\Sinv(\beta\up{3})\beta\up{2}}{(b\up{2} a S(b\up{3}))'}
    \\[-2pt]
    &\qquad\qquad\qquad\qquad\qquad\qquad\qquad\qquad\qquad
    \Smash
    (b\up{2} a S(b\up{3}))''b\up{4}\\
    &=
    \beta\up{3}(b\up{1}\leftreg\alpha)\Sinv(\beta\up{2})
    \beta\up{1}
    \Smash
    b\up{2} a S(b\up{3})b\up{4}
    \\
    &=
    \beta(b\up{1}\leftreg\alpha)
    \Smash
    b\up{2} a \ {}={} \ (\beta\Smash b)(\alpha\Smash a),
\end{align*}
where in $\smash[t]{{}\stackrel{\checkmark}{=}{}}$ we used that
$(a\rightreg\alpha)\leftreg\beta=\beta'\eval{\alpha\beta''}{a}$.

\subsection{Braided products} We now somewhat generalize the
observation leading to~\bref{thm:braided-comm}.  We first recall the
definition of a braided product, then see when braided commutativity
is hereditary under taking a braided product, and verify the
corresponding condition for $B^{*\mathrm{cop}}$ and $B$; their braided
product, which is therefore a braided commutative Yetter--Drin\-feld
module algebra, actually coincides with~$\HD(B^*)$.

\subsubsection{}
If $H$ is a Hopf algebra and $X$ and $Y$ two (left--left)
Yetter--Drin\-feld module algebras, their \textit{braided product}
$X\Braid Y$ is defined as the tensor product with the composition
\begin{equation}\label{braided-prod}
  (x\Braid y)(v\Braid u)
  = x(y\mone\leftii v)\Braid y\zero u,
  \quad x,v\in X,\quad y,u\in Y.
\end{equation}
This is a Yetter--Drin\-feld module algebra.  (Indeed, the
associativity of~\eqref{braided-prod} is ensured by $Y$ being a
comodule algebra and $X$ being a module algebra.  As a tensor product
of Yetter--Drin\-feld modules, $X\braid Y$ is a Yetter--Drin\-feld
module under the diagonal action (via iterated coproduct) and
codiagonal coaction of~$H$.  By the Yetter--Drin\-feld axiom for $Y$
and the module algebra properties of~$X$ and~$Y$, moreover, $X\braid
Y$ is a module algebra; the routine verification is given
in~\bref{standard} for completeness.  That $X\braid Y$ is a comodule
algebra follows from the comodule algebra properties of $X$ and $Y$
and the Yetter--Drin\-feld axiom for~$Y$; this is also recalled
in~\bref{standard}.)

\subsubsection{}\label{braided-sym}We say that two Yetter--Drin\-feld
modules $X$ and $Y$ are \textit{braided symmetric} if
\begin{equation*}
  c^{\vphantom{1}}_{Y,X}=c_{X,Y}^{-1}
\end{equation*}
(note that both sides here are maps $Y\tensor X\to X\tensor Y$), that
is,
\begin{equation*}
  (y\mone\leftii x)\tensor y\zero
  = x\zero\tensor\bigl(S^{-1}(x\mone)\leftii y\bigr).
\end{equation*}
\begin{lemma}\label{lemma:locked}
  Let $X$ and $Y$ be braided symmetric Yetter--Drin\-feld modules, each
  of which is a braided commutative Yetter--Drin\-feld module algebra.
  Then their braided product $X\Braid Y$ is also braided commutative.
\end{lemma}

We must show that
\begin{equation}\label{to-show}
  \bigl((x\Braid y)\mone\leftii(v\Braid u)\bigr)(x\Braid y)\zero
  =(x\Braid y)(v\Braid u)
\end{equation}
for all $x,v\in X$ and $y,u\in Y$.  For this, we write the condition
$c^{\vphantom{1}}_{X,Y}=c^{-1}_{Y,X}$ as
\begin{equation*}
  (x\mone\leftii y)\tensor x\zero
  = y\zero\tensor\bigl(S^{-1}(y\mone)\leftii x\bigr)
\end{equation*}
and use it to establish an auxiliary identity,
\begin{equation}\label{locked}
  \bigl((x\mone\leftii y)\mone\leftii x\zero\bigr)
  \tensor(x\mone\leftii y)\zero
  \begin{aligned}[t]
    &=\Bigl(y\zero{}\mone\leftii\bigl(S^{-1}(y\mone)\leftii x\bigr)
    \Bigr)\tensor y\zero{}\zero
    \\
    &=\bigl(y''\mone S^{-1}(y'\mone)\leftii x\bigr)\tensor y\zero
    =x\tensor y.
  \end{aligned}
\end{equation}
The left-hand side of~\eqref{to-show} is
\begin{multline*}
  \mbox{}\\*[-1.3\baselineskip]
  \shoveleft{
    \bigl((x\Braid y)\mone\leftii(v\Braid u)\bigr)(x\Braid y)\zero}
  \\
  \begin{aligned}[t]
    &=\bigl(x\mone y\mone\leftii(v\Braid u)\bigr)(x\zero\Braid y\zero)
    \\
    &=\bigl((x\mone' y\mone'\leftii v)\Braid
    (x\mone'' y\mone''\leftii u)\bigr)(x\zero\Braid y\zero)
    \\
    &=(x\mone' y\mone'\leftii v)
    \bigl((x\mone'' y\mone''\leftii u)\mone\leftii x\zero\bigr)
    \Braid
    (x\mone'' y\mone''\leftii u)\zero y\zero
    \\
    &=(x\mone y\mone'\leftii v)
    \bigl((x\zero{}\mone\leftii( y\mone''\leftii u))\mone\leftii
    x\zero{}\zero\bigr)
    \Braid
    (x\zero{}\mone\leftii(y\mone''\leftii u))\zero y\zero
    \\
    &=(x\mone y\mone'\leftii v)
    x\zero{}
    \Braid (y\mone''\leftii u)y\zero,
  \end{aligned}
\end{multline*}
just because of~\eqref{locked} in the last line.  But the right-hand
side of~\eqref{to-show} is
\begin{align*}
  (x\Braid y)(v\Braid u)
  &=x (y\mone\leftii v) \Braid y\zero u
  \\[-3pt]
  &=(x\mone y\mone\leftii v)x\zero
  \Braid (y\zero{}\mone\leftii u) y\zero{}\zero
\end{align*}
because $X$ and $Y$ are both braided commutative.  The two expressions
coincide.

\subsubsection{Remark}\label{opp-braid}Because the braided symmetry
condition is symmetric with respect to the two modules, we also have
the braided symmetric Yetter--Drin\-feld module algebra $Y\Braid X$,
with the product
\begin{equation*}
  (y\Braid x)(u\Braid v)=y(x\mone\leftii u)\Braid x\zero v.
\end{equation*}
In addition to the multiplication inside $Y$ and inside $X$, this
formula expresses the relations $x u = (x\mone\leftii u)x\zero$
satisfied in $Y\Braid X$ by $x\in X$ and $u\in Y$.  Because
$c^{\vphantom{1}}_{X,Y}=c_{Y,X}^{-1}$, these are the same relations
$ux = (u\mone\leftii x)u\zero$ that we have in $X\Braid Y$.  Somewhat
more formally, the isomorphism
\begin{equation*}
  \phi: X\Braid Y\to Y\Braid X
\end{equation*}
is given by $\phi:x\Braid y\mapsto (x\mone\leftii y)\Braid x\zero$.
This is a module map by virtue of the Yetter--Drin\-feld condition, and
it is immediate to verify that $\delta(\phi(x\Braid
y))=(\id\tensor\phi)(\delta(x\Braid y))$.
That 
$\phi$ is an algebra map follows by calculating
\begin{align*}
  \phi(x\Braid y)\phi(v\Braid u)
  &=((x\mone\leftii y)\Braid x\zero)
  ((v\mone\leftii u)\Braid v\zero)
  \\
  &=(x\mone\leftii y)(x\zero{}\mone v\mone\leftii u)
  \Braid x\zero{}\zero v\zero
  \\
  &=(x'\mone\leftii y)(x''\mone v\mone\leftii u)
  \Braid x\zero v\zero
  \\
  &=x\mone\leftii\bigl(y(v\mone\leftii u)\bigr)
  \Braid x\zero v\zero\\
  \intertext{and}
  \phi((x\Braid y)(v\Braid u))
  &=\phi\bigl(x(y\mone\leftii v)\Braid y\zero u\bigr)
  \\
  &=
  (x\mone(y\mone\leftii v)\mone\leftii(y\zero u))
  \Braid x\zero(y\mone\leftii v)\zero
  \\
  &\stackrel{\checkmark}{=}
  x\mone\leftii(y\zero u)\zero\Braid
  x\zero\bigl(S^{-1}(y\zero{}\mone u\mone)\leftii(y\mone\leftii v)\bigr)
  \\
  &=
  x\mone \leftii(y\zero u\zero)\Braid
  x\zero\bigl(S^{-1}(y''\mone u\mone)y'\mone\leftii v\bigr)
  \\
  &=
  x\mone\leftii(y u\zero)\Braid
  x\zero\bigl(S^{-1}(u\mone)\leftii v\bigr)
  \\
  &\stackrel{\checkmark}{=}
  x\mone\leftii\bigl(y(v\mone\leftii u)\bigr)
  \Braid
  x\zero v\zero,
\end{align*}
where the braided symmetry condition was used in each of the
${}\stackrel{\checkmark}{=}{}$ equalities.

\subsubsection{Multiple braided products}\label{XYXY} Further examples of
Yetter--Drinfeld module algebras are provided by multiple braided
products $X_{1}\Braid\dots\Braid X_{N}$ (of Yetter--Drin\-feld
$H$-module algebras $X_{i}$), defined as the corresponding tensor
products with the diagonal action and codiagonal coaction of $H$ and
with the relations
\begin{equation}\label{i>j}
  x[i]\Braid y[j] = (x\mone\leftii y)[j]\Braid x\zero[i]
  \ \ \text{for all}\ \ i>j,
\end{equation}
where $z[i]\in X_{i}$.  (The inverse relations are $x[i]\Braid y[j] =
y\zero[j]\Braid(S^{-1}(y\mone)\leftii x)[i]$, $i<j$.)  It readily
follows from the Yetter--Drin\-feld module algebra axioms for each of
the $X_i$ that $X_{1}\Braid\dots\Braid X_{N}$ is an associative
algebra and, in fact, a Yetter--Drin\-feld $H$-module algebra.

More specifically, let $X$ and $Y$ be braided symmetric
Yetter--Drin\-feld $H$-module algebras, as in~\bref{braided-sym}, and
consider the ``alternating'' products
\begin{equation}\label{eq:XYXY}
  X\Braid Y\Braid X\Braid Y\Braid\dots,
\end{equation}
with an arbitrary number of factors (or a similar product with the
leftmost $Y$, or actually their inductive limits).  We let $X[i]$
denote the $i$th copy of~$X$, and similarly with~$Y[j]$.  For
arbitrary $x[i]\in X[i]$ and $y[j]\in Y[j]$, we then have
the relations
\begin{equation}\label{allij}
  x[2i+1]\Braid y[2j] = (x\mone\leftii y)[2j]\Braid x\zero[2i+1],
\end{equation}
which by~\eqref{i>j} are satisfied for all $i\geq j$; but by the
braided symmetry condition, relations~\eqref{allij} \textit{hold for
  all $i$ and $j$} (replicating the relations between elements of $X$
and elements of $Y$ in~\hbox{$X\Braid Y$}).  In~\eqref{eq:XYXY}, also,
\begin{equation}\label{xxyy}
  \begin{aligned}
    x[2i+1]\Braid v[2j+1]&=(x\mone\leftii v)[2j+1]\Braid x\zero[2i+1],
    \quad x,v\in X,\\[-2pt]
    y[2i]\Braid u[2j]&=(y\mone\leftii u)[2j]\Braid y\zero[2i],
    \quad y,u\in Y,
  \end{aligned}\quad i > j.
\end{equation}
(These formulas also hold for $i=j$ if $X$ and $Y$ are braided
commutative.)

\subsection{$\HD(B^*)$ as a braided product}
Theorem~\bref{thm:braided-comm} can be reinterpreted by saying that
the Heisenberg double of $B^*$ is a braided product,
\begin{equation*}
  \HD(B^*) = B^{*\mathrm{cop}}\Braid B,
\end{equation*}
with the braiding
\begin{align*}
  b\tensor\beta&\mapsto (b\mone\leftii\beta)\tensor b\zero,
  \qquad b\in B,\quad\beta\in B^*,
\end{align*}
where we abbreviate the action of $B$ in~\bref{H-action} to
\begin{align*}
  m\leftii(\beta\Smash b)
  &=(m'\leftreg\beta)
  \Smash(m'' b S(m''')),\quad m\in B,
\end{align*}
and further use ${}\leftii{}$ for the restriction to $B^*$,
viz., $m\leftii\beta = m\leftreg\beta$.  It is also understood that
$B^{*\mathrm{cop}}$ and $B$ are viewed as left $\DD(B)$-comodule
algebras via
\begin{equation*}
  \delta:\beta\mapsto(\beta''\tensor 1)\tensor\beta',
  \qquad
  \delta: b\mapsto (\varepsilon\tensor b')\tensor b''
\end{equation*}
and left $\DD(B)$-module algebras via
\begin{equation*}
  (\mu\tensor m)\leftii\beta=\mu''(m\leftreg\beta)\Sinv(\mu'),
  \qquad
  (\mu\tensor m)\leftii b
  =(m' b S(m''))\rightreg\Sinv(\mu).
\end{equation*}
Both $B^{*\mathrm{cop}}$ and $B$ are then Yetter--Drin\-feld
$\DD(B)$-module algebras, and each is braided commutative.

Moreover, \textit{$B^{*\mathrm{cop}}$ and $B$ are braided symmetric}
because
$c^{\vphantom{1}}_{B^{*\mathrm{cop}},B}=c_{B,B^{*\mathrm{cop}}}^{-1}$,
i.e.,
\begin{equation*}
  (b\mone\leftii\beta)\tensor b\zero
  =\beta\zero\tensor(S_{_{\DD}}^{-1}(\beta\mone)\leftii b).
\end{equation*}
The antipode here is that of $\DD(B)$, and therefore the right-hand
side evaluates as $%
\beta'\tensor(S^{*}(\beta'')\leftii b) =\beta'\tensor(b\rightreg
\Sinv(S^{*}(\beta''))) =\beta'\tensor(b\rightreg \beta'')$, which is
immediately seen to coincide with the left-hand side.

Thus, the result that $\HD(B^*)=B^{*\mathrm{cop}}\Braid B$ is a
braided commutative Yetter--Drin\-feld module algebra now follows
from~\bref{lemma:locked}.  (This offers a nice alternative to an
unilluminating brute-force proof.)

\subsection{Heisenberg $n$-tuples$/$chains}\label{sec:chains}
It follows from~\bref{opp-braid} that $\HD(B^*)$ is also isomorphic to
the braided commutative Yetter--Drin\-feld module algebra $B\Braid
B^{*\mathrm{cop}}$, with the product
\begin{equation*}
  (a\Braid\alpha)(b\Braid\beta)=a(b\rightreg\Sinv(\alpha''))
  \Braid\alpha'\beta.
\end{equation*}

We next consider ``Heisenberg $n$-tuples$/$chains''---\,the
alternating products
\begin{align*}
  \HH_{2n} &=
  B^{*\mathrm{cop}}\Braid B\Braid B^{*\mathrm{cop}}\Braid B
  \Braid\dots\Braid B,\\[-2pt]
  \HH_{2n+1} &=
  B^{*\mathrm{cop}}\Braid B\Braid B^{*\mathrm{cop}}\Braid B
  \Braid \dots \Braid B\Braid B^{*\mathrm{cop}}.
\end{align*}
As we saw in~\bref{XYXY}, the relations are then given by
\begin{align*}
  b[2i]\,\beta[2j+1] &= (b'\leftreg\beta)[2j+1]\,b''[2i]
  \quad\text{for all $i$ and $j$}
\end{align*}
(where $b\in B$ and $\beta\in B^*$, $B^{*\mathrm{cop}}\to
B^{*\mathrm{cop}}[2j+1]$ and $B\to B[2i]$\pagebreak[3] are the
morphisms onto the respective factors, and we omit ${}\Braid{}$ for
brevity)
and
\begin{equation*}
  \begin{aligned}
    \alpha[2i+1]\,\beta[2j+1]&=(\alpha'''\beta\Sinv(\alpha''))[2j+1]\,
    \alpha'[2i+1],
    \quad \alpha,\beta\in B^{*\mathrm{cop}},\\[-2pt]
    a[2i]\, b[2j]&=(a' b S(a''))[2j]\, a'''[2i],
    \quad a,b\in B,
  \end{aligned}\quad i\geq j
\end{equation*}
(Relations inverse to the last two are
$\beta[2i+1]\,\alpha[2j+1]=\alpha'[2j+1]
(S^*(\alpha'')\beta\alpha''')[2i+1]$ and $b[2i]\,
a[2j]=a'''[2j]\,(S^{-1}(a'')b a')[2i]$ for $i\leq j$.)

The chains with the leftmost $B$ factor are defined entirely
similarly.  The obvious embeddings allow defining (one-sided or
two-sided) inductive limits of alternating chains.  All the chains are
Yetter--Drin\-feld module algebras, but those with $\geq3$ tensor
factors are not braided commutative in general.

\section{Yetter--Drinfeld module algebras for $\UresSL2$}
\label{sec:sl2}
In this section, we construct Yetter--Drin\-feld module algebras for
$\UresSL2$ at an even root of unity
\begin{equation*}
  \q=e^{\frac{i\pi}{p}}
\end{equation*}
for an integer $p\geq2$.  $\UresSL2$ is the $2p^3$-dimensional quantum
group with generators $E$, $K$, and $F$ and the relations
\begin{gather*}
  KEK^{-1}=\q^2E,\quad
  KFK^{-1}=\q^{-2}F,\quad
  [E,F]=\ffrac{K-K^{-1}}{\q-\q^{-1}},
  \\
  E^{p}=F^{p}=0,\quad K^{2p}=1
\end{gather*}
and the Hopf algebra structure $\Delta(E)= E\tensor K + 1\tensor E$,
$\Delta(K)=K\tensor K$, $\Delta(F)=F\tensor 1 + K^{-1}\tensor F$,
$\epsilon(E)=\epsilon(F)=0$, $\epsilon(K)=1$, $S(E)=-E K^{-1}$,
$S(K)=K^{-1}$, $S(F)=-K F$.\footnote{In an ``applied'' context (see,
  e.g., \cite{[S-q],[AM],[NT]}), this quantum group first appeared
  in~\cite{[FGST],[FGST2]}; subsequently, it gradually transpired
  (with the final picture having emerged from~\cite{[KS]}) that that
  was just a continuation of a series of previous
  (re)discoveries~\cite{[AGL],[Su],[X]} (also see~\cite{[Erd]}).  The
  ribbon and (somewhat stretching the definition) factorizable
  structures of $\UresSL2$ were worked out in~\cite{[FGST]}.}
 
In~\cite{[FGST],[FGST2]}, $\UresSL2$ was arrived at as a subquotient
of the Drinfeld double of a Taft Hopf algebra (a trick also used,
e.g., in~\cite{[Sch-Galois]} for a closely related quantum group).  It
turns out that not only $\DD(B)$ but also the pair $(\DD(B),\HD(B^*))$
can be ``truncated'' to a pair $(\UresSL2, \HresSL2)$ of
$2p^3$-dimensional algebras, with $\HresSL2$ being a braided
commutative Yetter--Drin\-feld $\UresSL2$-module algebra.  This is
worked out in what follows.  $\HresSL2$\,---\,a ``Heisenberg
counterpart'' of $\UresSL2$\,---\,appears in~\bref{sec:Hq}.

\subsection{$\DD(B)$ and $\HD(B^*)$ for the Taft Hopf algebra $B$}
\subsubsection{The Taft Hopf algebra $B$}
Let 
\begin{equation*}
  B=\mathrm{Span}(E^m k^n),\quad
  0\leq m\leq p-1,\quad 0\leq n\leq 4p-1,
\end{equation*}
be the $4p^2$-dimensional Hopf algebra generated by~$E$ and~$k$ with
the relations\pagebreak[3]%
\begin{gather}\label{prod-B}
  k E =\q E k,\quad E^p=0,\quad k^{4p}=1,
\end{gather}
and with the comultiplication, counit, and antipode given by
\begin{gather}\label{coalgebra-B}
  \begin{gathered}
    \Delta(E)=1\otimes E+ E\otimes k^2,\quad
    \Delta(k)= k\otimes k,\\
    \epsilon(E)=0,\quad\epsilon(k)=1,\\
    S(E)=- E k^{-2},\quad S(k)= k^{-1}.
  \end{gathered}
\end{gather}

Dual elements $F,\varkappa\in B^*$ are introduced as
\begin{equation*}
  \eval{F}{E^{m} k^{n}}=\delta_{m,1}\ffrac{\q^{-n}}{\q-\q^{-1}},
  \qquad
  \eval{\varkappa}{E^{m} k^{n}}=\delta_{m,0}\q^{-n/2}.
\end{equation*}
Then~\cite{[FGST]}
\begin{equation*}
  B^*=\mathrm{Span}(F^a\varkappa^b),\quad
  0\leq a\leq p-1,\quad 0\leq b\leq 4p-1.
\end{equation*}

\subsubsection{The Drinfeld double $\DD(B)$}
Straightforward calculation shows~\cite{[FGST]} that the Drinfeld
double $\DD(B)$ is the Hopf algebra generated by $E$, $F$, $k$, and
$\varkappa$ with the relations given~by
\begin{itemize}
\item[i)] relations~\eqref{prod-B} in~$B$,
\item[ii)] the relations
  \begin{equation*}
    \varkappa F=\q F\varkappa,\quad F^p=0,\quad
    \varkappa^{4p}=1
  \end{equation*}
  in $B^*$, and
\item[iii)] the cross-relations
  \begin{gather*}
    k\varkappa=\varkappa k,\quad k F k^{-1}=\q^{-1} F,\quad
    \varkappa E\varkappa^{-1}=\q^{-1} E,\quad
    [E, F]=\mfrac{ k^2-\varkappa^2}{\q-\q^{-1}}.
  \end{gather*}
\end{itemize}
The Hopf-algebra structure
$(\Delta_{_{\DD}},\varepsilon_{_{\DD}},S_{_{\DD}})$ of $\DD(B)$ is
given by~\eqref{coalgebra-B} and
\begin{gather*}
  \Delta_{_{\DD}}(F)=\varkappa^2\tensor F+ F\tensor1,\quad
  \Delta_{_{\DD}}(\varkappa)=\varkappa\tensor\varkappa,\quad
  \epsilon_{_{\DD}}(F)=0,\quad
  \epsilon_{_{\DD}}(\varkappa)=1,
  \\
  S_{_{\DD}}(F)=-\varkappa^{-2} F,
  \quad S_{_{\DD}}(\varkappa)=\varkappa^{-1}.
\end{gather*}

\subsubsection{The Heisenberg double $\HD(B^*)$}\label{HD-relations}
For the above $B$, \ $\HD(B^*)$ is spanned by
\begin{equation}\label{HBstar-span}
  F^a\varkappa^b\Smash E^c k^d,\qquad a,c=0,\dots,p-1,\quad
  b,d\in\oZ/(4p\oZ),
\end{equation}
where $\varkappa^{4p}=1$, $k^{4p}=1$, $F^p=0$, and $E^p=0$.  Then the
product in~\eqref{H-comp} becomes~\cite{[S-H]}
\begin{multline}\label{prod-in-smash-full}
  (F^r \varkappa^s\Smash E^m k^n)(F^a\varkappa^b\Smash E^c k^d)
  \\
  {}= \sum_{u\geq 0}\q^{-\half u(u-1)}\qbin{m}{u}
  \qbin{a}{u}\ffrac{\qfac{u}}{(\q - \q^{-1})^u}\,
  \q^{-\half b n + c n + a (s - n) +  u(2 c - a - b  + m - s)}
  \\*[-6pt]
  {}\times F^{a+r-u}\varkappa^{b+s} \Smash E^{m + c - u}k^{n+d+2u}.
\end{multline}

A convenient basis in $\HD(B^*)$ can be chosen as $(\varkappa, z,
\lambda, \Dz)$, where $\varkappa$ is understood as $\varkappa\Smash1$
and
\begin{align*}
  z &= -(\q-\q^{-1}) \varepsilon\Smash E k^{-2},\\
  \lambda &=\varkappa \Smash k,\\
  \Dz &= (\q-\q^{-1}) F\Smash 1.
\end{align*} 
The relations in $\HD(B^*)$ then become
  $\varkappa z = \q^{-1} z \varkappa$,
  $\varkappa \lambda = \q^{\half} \lambda \varkappa$,
  $\varkappa \Dz = \q \Dz \varkappa$,
    $\varkappa^{4p}=1$, and
\begin{gather*}
 \lambda^{4p}=1,\qquad
  z^p=0,\qquad \Dz^p=0,
  \\
  \lambda z = z \lambda,\qquad   \lambda \Dz = \Dz \lambda,\\
  \Dz z = (\q-\q^{-1}) 1 +  \q^{-2} z\Dz.
\end{gather*}

\subsection{The $(\UresSL2,\HresSL2)$ pair}
\subsubsection{From $\DD(B)$ to $\UresSL2$}
The ``truncation'' whereby $\DD(B)$ yields $\UresSL2$~\cite{[FGST]}
consists of two steps: first, taking the quotient
\begin{gather}\label{quotient}
  \bDDB=\DD(B)/(\varkappa k - 1)
\end{gather}
by the Hopf ideal generated by the central element $\varkappa\tensor k
- \varepsilon\tensor 1$ and, second, identifying $\UresSL2$ as the
subalgebra in $\bDDB$ spanned by $F^{\ell} E^{m} k^{2n}$ (tensor
product omitted) with $\ell,m=0,\dots,p-1$ and $n=0,\dots,2p-1$.  It
then follows that $\UresSL2$ is a Hopf algebra\,---\,the one described
at the beginning of this section, where $K=k^2$.

\subsubsection{From $\HD(B^*)$ to $\HresSL2$}\label{sec:Hq}
In $\HD(B^*)$, dually, we take a subalgebra and then a
quotient~\cite{[S-H]}.  In the basis chosen above, the subalgebra
(which is also a $\UresSL2$ submodule) is the one generated by $z$,
$\Dz$, and $\lambda$.  Its quotient by $\lambda^{2p}=1$ gives a
\textit{$2p^3$-dimensional algebra $\HresSL2$\,---\,the ``Heisenberg
  counterpart'' of \,$\UresSL2$}~\cite{[S-H]}.

As an associative algebra,
\begin{equation*}
  \HresSL2 = \Czd\tensor(\oC[\lambda]/(\lambda^{2p}-1)),
\end{equation*}
with the $p^2$-dimensional algebra
\begin{equation*}
  \Czd=\oC[z,\Dz]/(z^p, \ \Dz^p, \ \Dz z - (\q-\q^{-1}) - \q^{-2} z\Dz).
\end{equation*}

The $\UresSL2$ action on $\HresSL2$ follows from~\eqref{the-action} as
\begin{alignat*}{3}
  E\leftii\lambda^n&=\q^{-\frac{n}{2}}\qint{\fffrac{n}{2}}\,\lambda^n\,z,
  &
  k^2\leftii\lambda^n&=\q^{-n}\lambda,
  &
  F\leftii\lambda^n&=-\q^{\frac{n}{2}}\qint{\fffrac{n}{2}}\,\lambda^n\,\Dz,
  \quad
  \\
  E\leftii z^n &=-\q^n \qint{n} z^{n+1},
  &
  k^2\leftii z^n &= \q^{2n}\,z^n,
  &
  F\leftii z^n &= \qint{n} \q^{1-n}\,z^{n-1},
  \\
  E\leftii \Dz^n &= \q^{1-n}\qint{n}\Dz^{n-1},\quad
  &
  k^2\leftii \Dz^n &= \q^{-2 n} \Dz^n,\quad
  &
  F\leftii \Dz^n &= -\q^n \qint{n} \Dz^{n+1}.
\end{alignat*}

The coaction $\delta:\HresSL2\to\UresSL2\tensor\HresSL2$ follows
from~\eqref{delta} as
\begin{align*}
  \lambda&\mapsto1\tensor\lambda,\\
  z^m &\mapsto \sum_{s=0}^m(-1)^s\q^{s(1-m)}(\q - \q^{-1})^s
  \qbin{m}{s}\,
  E^s k^{-2m}\tensor z^{m-s},\\
  \Dz^m &\mapsto \sum_{s=0}^m\q^{s(m-s)}(\q - \q^{-1})^s
  \qbin{m}{s}\,
  F^s k^{-2(m-s)}\tensor \Dz^{m-s}.
\end{align*}
In particular,
$z \mapsto k^{-2}\tensor z -(\q - \q^{-1}) E k^{-2}\tensor 1$ and $\Dz
\mapsto k^{-2}\tensor \Dz +(\q - \q^{-1}) F \tensor 1$.

\subsubsection{} \textit{With the $\UresSL2$ action and coaction given
  above, $\HresSL2$ is a braided commutative Yetter--Drin\-feld
  $\,\UresSL2$-module algebra}.

Hence, in particular, $\Czd$ is also a braided commutative
Yetter--Drin\-feld $\,\UresSL2$-module algebra.\footnote{We also
  recall that $\Czd$ is in fact $\Mat_p(\oC)$~\cite{[S-U]}.}

\subsection{Heisenberg ``chains.''}
The Heisenberg $n$-tuples$/$chains defined in~\bref{sec:chains} can
also be ``truncated'' similarly to how we passed from $\HD(B^*)$ to
$\HresSL2$.  An additional possibility here is to drop the coinvariant
$\lambda$ altogether, which leaves us with the ``\textit{truly
  Heisenberg}'' Yetter--Drin\-feld $\UresSL2$-module algebras
\begin{align*}
  \HHH_2 &= \oC_{\q}^{*p}[\Dz_1]\Braid\oC_{\q}^p[z_2]
  = \oC_{\q}[z_2,\Dz_1],\\
  \HHH_{2n} &= \oC_{\q}^{*p}[\Dz_1]\Braid\oC_{\q}^p[z_2]\Braid\dots
  \Braid\oC_{\q}^{*p}[\Dz_{2n-1}]\Braid\oC_{\q}^p[z_{2n}],\\
  \HHH_{2n+1} &= \oC_{\q}^{*p}[\Dz_1]\Braid\oC_{\q}^p[z_2]\Braid\dots
  \Braid\oC_{\q}^{*p}[\Dz_{2n-1}]\Braid\oC_{\q}^p[z_{2n}]
  \Braid\oC_{\q}^{*p}[\Dz_{2n+1}]
\end{align*}
(or their infinite versions), where
$\oC_{\q}^{*p}[\Dz]=\oC[\Dz]/\Dz^p$ and $\oC_{\q}^{p}[z]=\oC[z]/z^p$,
with the braiding inherited from~\bref{sec:chains}, which amounts to
using the relations
\begin{align*}
  &\Dz_i\, z_{j} = \q - \q^{-1} + \q^{-2} z_{j}\,\Dz_i
  \\[-4pt]
  \intertext{for \textit{all} (odd) $i$ and (even) $j$, and}
  &\begin{aligned}
    z_i\,z_j &= \smash[t]{\q^{-2} z_j\,z_i + (1 - \q^{-2}) z_j^2},\\
    \Dz_i\,\Dz_j &= \q^2 \Dz_j\,\Dz_i + (1 - \q^2)\Dz_j^2,
  \end{aligned}\quad i\geq j
\end{align*}
(and $z_i^p=0$ and $\Dz_i^p=0$; our relations may be interestingly
compared with those in para-Grassmann algebras in~\cite{[FIK1]}).

\subsubsection*{Acknowledgments} I am grateful to A.~Isaev for the
useful comments.  This work was supported in part by the RFBR grant
07-01-00523, the RFBR--CNRS grant 09-01-93105, and the grant
LSS-1615.2008.2.

\appendix
\section{Drinfeld double}\label{app:D-double}
We recall that the Drinfeld double of $B$, denoted by $\DD(B)$, is
$B^*\tensor B$ as a vector space, endowed with the structure of a
quasitriangular Hopf algebra given as follows.  The co\-algebra
structure is that of $B^{*\mathrm{cop}}\tensor B$, the algebra
structure is given by
\begin{equation}\label{in-double}
  (\mu\tensor m)(\nu\tensor n)
  = \mu(m'\leftreg \nu\rightreg S^{-1}(m'''))
  \tensor m''n
\end{equation}
for all $\mu,\nu\in B^*$ and $m,n\in B$, the antipode is given by
\begin{equation}\label{antipode-double}
  S_{_{\DD}}(\mu\tensor m)
  =(\varepsilon\tensor S(m))(\Sinv(\mu)\tensor1)
  =(S(m''')\leftreg\Sinv(\mu)\rightreg m')\tensor S(m''),
\end{equation}
and the universal $R$-matrix is
\begin{equation}\label{R-double}
  R= \sum_I(\varepsilon \tensor e_I)\tensor (e^I\tensor 1),
\end{equation}
where $\{e_I\}$ is a basis of $B$ and $\{e^I\}$ its dual basis
in~$B^*$.

\section{Standard calculations}
\subsection{Proof of the action
  in~\eqref{the-action}}\label{app-prove-Drinfeld}
To show that~\eqref{the-action} defines an action of $\DD(B)$, we
verify~\eqref{to-show-action} by first evaluating its right-hand side:
\begin{align*}
  \bigl((m'\leftreg\mu\rightreg S^{-1}(m'''))\tensor m''\bigr)
    \leftii(\alpha\Smash a)\kern-175pt
  \\
  &=
  \eval{\mu'}{S^{-1}(m\up{5})}
  \eval{\mu'''}{m\up{1}}(\mu'' \tensor 1)
  \leftii\bigl((m\up{2}\leftreg\alpha)\Smash m\up{3} a
  S(m\up{4})\bigr)\\
  &=
  \eval{\mu\up{1}\!}{S^{-1}(m\up{5})}
  \eval{\mu\up{5}\!}{m\up{1}}
  \!\mu\up{4}(m\up{2}\!\leftreg\alpha)\Sinv(\mu\up{3})\Smash
  (m\up{3} a S(m\up{4})\rightreg\Sinv(\mu\up{2}))\kern-40pt
  \\
  &=
  \eval{\mu\up{1}}{S^{-1}(m\up{7})}
  \eval{\mu\up{5}\!}{m\up{1}}
  \mu\up{4}(m\up{2}\leftreg\alpha)\Sinv(\mu\up{3})\Smash m\up{4} a''
  S(m\up{5})\\[-3pt]
  &\qquad\qquad\qquad\qquad\qquad\qquad\qquad\qquad\qquad\qquad\quad
  {}\times\eval{\mu\up{2}}{m\up{6} S^{-1}(a')S^{-1}(m\up{3})}\\
  &=
  \eval{\mu\up{1}\!}{S^{-1}(a')S^{-1}(m\up{3})}
  \eval{\mu\up{4}\!}{m\up{1}}
  \mu\up{3}(m\up{2}\leftreg\alpha)\Sinv(\mu\up{2})\Smash
  m\up{4} a'' S(m\up{5})\\
  &=
  \eval{\Sinv(\mu\up{1})}{a'}
  \eval{\Sinv(\mu\up{2})}{m\up{3}}
  \eval{\mu\up{5}\!}{m\up{1}}\eval{\alpha''\!}{m\up{2}}
  \mu\up{4}\alpha'\Sinv(\mu\up{3})\\*[-3pt]
  &\qquad\qquad\qquad\qquad\qquad\qquad\qquad\qquad\qquad\qquad\qquad
  \qquad\qquad\quad
  \Smash m\up{4} a'' S(m\up{5})\\
  &=
  \eval{\Sinv(\mu\up{1})}{a'}
  \eval{\mu\up{5}\alpha''\Sinv(\mu\up{2})}{m\up{1}}
  \mu\up{4}\alpha'\Sinv(\mu\up{3})\Smash m\up{2} a'' S(m\up{3})\\
  &=
  \eval{\Sinv(\mu\up{1})}{a'}
  \eval{\bigl(\mu\up{3}\alpha\Sinv(\mu\up{2})\bigr)''}{m\up{1}}
  \bigl(\mu\up{3}\alpha\Sinv(\mu\up{2})\bigr)'\Smash m\up{2} a''
  S(m\up{3})
  \\
  &=
  \bigl(m\up{1}\leftreg(\mu\up{3}\alpha\Sinv(\mu\up{2}))\bigr)\Smash
  m\up{2}\bigl(a\rightreg \Sinv(\mu\up{1})\bigr)S(m\up{3}),
\end{align*}%
which is the same as the left-hand side:
\begin{align*}
  (\varepsilon\tensor m)\leftii
  \bigl((\mu\tensor 1)\leftii(\alpha\Smash a)\bigr)
  &=
  (\varepsilon\tensor m)\leftii
  \bigl(\mu'''\alpha\Sinv(\mu'')
  \Smash(a\rightreg\Sinv(\mu'))\bigr)\\
  &=\bigl(m'\leftreg(\mu'''\alpha\Sinv(\mu''))\bigr)
  \Smash\bigl(m''(a\rightreg\Sinv(\mu'))S(m''')\bigr).
\end{align*}

\subsection{Proof of the $\DD(B)$-module algebra property for $\HD(B^*)$}
\label{app-prove-module-algebra}
To show \eqref{m-a} for the action in~\eqref{the-action}, we do this
for $M=\varepsilon\tensor m$ and $M=\mu\tensor 1$ separately.

First, the right-hand side of~\eqref{m-a} with $M=\varepsilon\tensor
m$ is
\begin{multline*}
  \mbox{}\\[-1.2\baselineskip]
  \shoveleft{\bigl((\varepsilon\tensor m')\leftii(\alpha\Smash a)\bigr)
    \bigl((\varepsilon\tensor m'')\leftii(\beta\Smash b)\bigr)}\\
  \begin{aligned}[t]
    &=
    \bigl((m\up{1}\leftreg\alpha)\Smash m\up{2} a S(m\up{3})\bigr)
    \bigl((m\up{4}\leftreg\beta)\Smash m\up{5} b S(m\up{6})\bigr)
    \\
    &=
    (m\up{1}\leftreg\alpha)
    \bigl(((m\up{2} a S(m\up{3}))'m\up{4})\leftreg\beta\bigr)
    \Smash
    (m\up{2} a S(m\up{3}))''m\up{5} b S(m\up{6})
    \\
    &=
    (m\up{1}\leftreg\alpha)
    \bigl(m\up{2} a'\leftreg\beta\bigr)
    \Smash
    m\up{3} a'' S(m\up{4}) m\up{5} b S(m\up{6})
    \\
    &=
    (m\up{1}\leftreg\alpha)
    \bigl(m\up{2} a'\leftreg\beta\bigr)
    \Smash
    m\up{3} a''  b S(m\up{4}),
  \end{aligned}
\end{multline*}
which is the left-hand side $(\varepsilon\tensor m)
\leftii\bigl(\alpha(a'\leftreg\beta) \Smash a''b\bigr)$.

Second, the left-hand side of~\eqref{m-a} with $M=\mu \tensor 1$ is
\begin{multline*}
  \mbox{}\\[-1.2\baselineskip]
  \shoveleft{
    (\mu\tensor 1)\leftii\bigl(\alpha(a'\leftreg\beta)\Smash
    a''b\bigr)}
  \\
  \begin{aligned}[t]
    &=\mu'''\alpha(a'\leftreg\beta)\Sinv(\mu'')\Smash
    \bigl((a''b)\rightreg\Sinv(\mu')\bigr)\\
    &=\mu\up{4}\alpha(a'\leftreg\beta)\Sinv(\mu\up{3})\Smash
    (a''\rightreg\Sinv(\mu\up{2}))(b\rightreg\Sinv(\mu\up{1})).
  \end{aligned}
\end{multline*}
But the right-hand side of~\eqref{m-a} evaluates the same:
\begin{align*}
  \bigl((\mu''\tensor 1)\leftii(\alpha\Smash a)\bigr)
  \bigl((\mu'\tensor 1)\leftii(\beta\Smash b)\bigr)\kern-160pt
  \\
    &=
    \bigl(\mu\up{6}\alpha\Sinv(\mu\up{5})
    \Smash
    (a\rightreg\Sinv(\mu\up{4}))\bigr)
    \bigl(\mu\up{3}\beta\Sinv(\mu\up{2})
    \Smash
    (b\rightreg\Sinv(\mu\up{1}))\bigr)
    \\
    &=
    \mu\up{6}\alpha\Sinv(\mu\up{5})
    \bigl((a\rightreg\Sinv(\mu\up{4}))'\leftreg
    \mu\up{3}\beta\Sinv(\mu\up{2})\bigr)
    \\[-3pt]
    &\qquad\qquad\qquad\qquad\qquad\qquad\qquad\qquad
    {}\Smash
    \bigl(a\rightreg\Sinv(\mu\up{4})\bigr)''
    \bigl(b\rightreg\Sinv(\mu\up{1})\bigr)
    \\
    &=
    \mu\up{6}\alpha\Sinv(\mu\up{5})
    \bigl((a'\rightreg\Sinv(\mu\up{4}))\leftreg
    (\mu\up{3}\beta\Sinv(\mu\up{2}))\bigr)
    \\[-2pt]
    &\qquad\qquad\qquad\qquad\qquad\qquad\qquad\qquad
    {}\Smash a'' \bigl(b\rightreg\Sinv(\mu\up{1})\bigr)\\*[-3pt]
    &\phantom{{}={}}
    \text{\mbox{}\hfill(because $\Delta(a\rightreg\mu)
      =(a'\rightreg\mu)\tensor a''$)}
    \\[1pt]
    &=
    \mu\up{6}\alpha\Sinv(\mu\up{5})
    \eval{\Sinv(\mu\up{4})(\mu\up{3}\beta\Sinv(\mu\up{2}))''}{a'}
    (\mu\up{3}\beta\Sinv(\mu\up{2}))'
    \\[-3pt]
    &\qquad\qquad\qquad\qquad\qquad\qquad\qquad\qquad
    {}\Smash
    a'' \bigl(b\rightreg\Sinv(\mu\up{1})\bigr)\\*[-3pt]
    &\phantom{{}={}}
    \text{(simply because $(a\rightreg\alpha)\leftreg\beta
      =\beta'\eval{\alpha\beta''}{a}$)}
    \\[1pt]
    &=
    \eval{\beta''\Sinv(\mu\up{2})}{a'}
    \mu\up{4}\alpha\beta'\Sinv(\mu\up{3})
    {}\Smash
    a'' \bigl(b\rightreg\Sinv(\mu\up{1})\bigr)
    \\
    &=
    \eval{\beta''}{a'}
    \eval{\Sinv(\mu\up{2})}{a''}
    \mu\up{4}\alpha\beta'\Sinv(\mu\up{3})
    {}\Smash
    a''' \bigl(b\rightreg\Sinv(\mu\up{1})\bigr)
    \\
    &=
    \mu\up{4}\alpha(a'\leftreg\beta)\Sinv(\mu\up{3})
    {}\Smash
    (a''\rightreg \Sinv(\mu\up{2}))
    \bigl(b\rightreg\Sinv(\mu\up{1})\bigr).
\end{align*}

\subsection{Standard checks for braided products}\label{standard}
Here, we give the standard calculations establishing the module
algebra and comodule algebra properties for the product defined
in~\eqref{braided-prod}.  

The module algebra property follows by calculating
\begin{align*}
  \bigl(h'\leftii (x\Braid y)\bigr)
  \bigl(h''\leftii (v\Braid u)\bigr)
    & =\bigl((h\up{1}\leftii x)\Braid(h\up{2}\leftii y)\bigr)
    \bigl((h\up{3}\leftii v)\Braid (h\up{4}\leftii u)\bigr)
    \\
    & =(h\up{1}\leftii x)
    \bigl((h\up{2}\leftii y)\mone h\up{3}\leftii v\bigr)
    \Braid (h\up{2}\leftii y)\zero (h\up{4}\leftii u)
    \\
    & =(h\up{1}\leftii x)
    (h\up{2} y\mone \leftii v)
    \Braid (h\up{3}\leftii y\zero) (h\up{4}\leftii u)
    \\
    & = 
    h\leftii\bigl(x(y\mone\leftii v)\Braid y\zero u\bigr)
    =h\leftii\bigl((x\Braid y) (v\Braid u)\bigr).
\end{align*}

To verify the comodule algebra property $\delta\bigl((x\Braid
y)(v\Braid u)\bigr)= \delta(x\Braid y)\delta(v\Braid u)$, we
calculate the left-hand side using that $X$ and $Y$ are comodule
algebras and that $Y$ is Yetter--Drin\-feld:
\begin{align*}
  \delta\bigl((x\Braid y)(v\Braid u)\bigr)
  &= \delta\bigl(x(y\mone\leftii v)\Braid y\zero u\bigr)
  \\
  &=\bigl(x(y\mone\leftii v)\bigr)\mone(y\zero u)\mone
  \tensor
  \bigl(x(y\mone\leftii v)\bigr)\zero\Braid\bigl(y\zero u\bigr)\zero
  \\
  &=x\mone(y\mone\leftii v)\mone y\zero{}\mone u\mone
  \tensor
  x\zero(y\mone\leftii v)\zero
  \Braid y\zero{}\zero u\zero
  \\
  &=x\mone(y\mone'\leftii v)\mone y\mone'' u\mone
  \tensor
  \bigl(x\zero(y\mone'\leftii v)\zero
  \Braid y\zero u\zero\bigr)
  \\
  &=x\mone y\mone' v\mone u\mone
  \tensor
  \bigl(x\zero (y\mone''\leftii v\zero)
  \Braid y\zero u\zero\bigr),\\[-3pt]
  \intertext{which is the same as the right-hand side by another use
    of the comodule axiom for~$Y$:}
  \delta(x\Braid y)\delta(v\Braid u)
  &=\bigl(x\mone y\mone\tensor(x\zero\Braid y\zero)\bigr)
  \bigl(v\mone u\mone\tensor(v\zero\Braid u\zero)\bigr)
  \\
  &=(x\mone y\mone v\mone u\mone)
  \tensor
  \bigl(x\zero (y\zero{}\mone\leftii v\zero)\Braid y\zero{}\zero u\zero
  \bigr)
  \\
  &=(x\mone y\mone' v\mone u\mone)
  \tensor
  \bigl(x\zero (y\mone''\leftii v\zero)\Braid y\zero u\zero
  \bigr).
\end{align*}

\parindent=0pt


\begin{thebibliography}{99}
\bibitem{[AF]} A.Yu.\;Alekseev and L.D.\;Faddeev, \textit{$(T^*G)_t$:
    A toy model for conformal field theory}, Commun.\ Math.\ Phys. 141
  (1991) 413--422.


\bibitem{[RSts]} N.Yu.\;Reshetikhin, and M.A.\;Semenov-Tian-Shansky,
  \textit{Central extensions of quantum current groups}, Lett.\ Math.\
  Phys. 19 (1990) 133--142.


\bibitem{[Sts]}M.A.\;Semenov-Tyan-Shanskii, \textit{Poisson--Lie
    groups. The quantum duality principle and the twisted quantum
    double},
  Theor.\ Math.\ Phys.\ 93 (1992) 1292--1307.


\bibitem{[Lu-double]}J.-H.\;Lu, \textit{On the Drinfeld double and the
    Heisenberg double of a Hopf algebra}, Duke Math. J. 74 (1994)
  763--776.


\bibitem{[Lu-alg]}J.-H.\;Lu, \textit{Hopf algebroids and quantum
    groupoids}, math.QA$/$\linebreak[0]9505024.


\bibitem{[P]}F.\;Panaite, \textit{Doubles of (quasi) Hopf algebras and
    some examples of quantum groupoids and vertex groups related to
    them}, math.QA$/$\linebreak[0]0101039.


\bibitem{[BM]} T.\;Brzezi\'nski and G.\;Militaru,
  \textit{Bialgebroids, $\times_{A}$-bialgebras and duality}, J.\
  Algebra 251 (2002) 279--294 [math.QA$/$\linebreak[0]0012164].


\bibitem{[K]} R.M.\;Kashaev, \textit{Heisenberg double and the
    pentagon relation},
  St.\ Petersburg Math.\ J. 8 (1997) 585--592
  [q-alg$/$\linebreak[0]9503005].


\bibitem{[Ka]}M.\;Kapranov, \textit{Heisenberg doubles and derived
    categories}, q-alg$/$\linebreak[0]9701009.


\bibitem{[Mi]}G.\;Militaru, \textit{Heisenberg double, pentagon
    equation, structure and classification of finite-dimensional Hopf
    algebras}, J.\ London Math. Soc. (2) 69 (2004) 44--64.


\bibitem{[IP]}A.\;Isaev and P.\;Pyatov, \textit{Spectral extension of
    the quantum group cotangent bundle}, Comm.\ Math.\ Phys.\ 288
  (2009) 1137--1179 [arXiv:\linebreak[0]0812.2225 [math.QA]].


\bibitem{[FGST]} B.L.\;Feigin, A.M.\;Gainutdinov, A.M.\;Semikhatov,
  and I.Yu.\;Tipunin, \textit{Modular group representations and fusion
    in logarithmic conformal field theories and in the quantum group
    center}, Commun.\ Math.\ Phys.\ 265 (2006) 47--93
  [hep-th$/$\linebreak[0]0504093].


\bibitem{[FGST2]} B.L.\;Feigin, A.M.\;Gainutdinov, A.M.\;Semikhatov,
  and I.Yu.\;Tipunin, \textit{Kazhdan--Lusztig correspondence for the
    representation category of the triplet $W$-algebra in logarithmic
    CFT}, Theor.\ Math.\ Phys.\ 148 (2006) 1210--1235
  [math.QA$/$\linebreak[0]0512621].


\bibitem{[S-q]} A.M.\;Semikhatov, \textit{Factorizable ribbon quantum
    groups in logarithmic conformal field theories}, Theor.\
  Math. Phys.\ 154 (2008) 433--453 [arXiv:\linebreak[0]0705.4267
  [hep-th]].


\bibitem{[G]}A.M.\;Gainutdinov, \textit{A generalization of the
    Verlinde formula in logarithmic conformal field theory},
  Theor.\ Math. Phys. 159 (2009) 575--586.


\bibitem{[S-U]}A.M.\;Semikhatov, \textit{A differential
    $\mathscr{U}$-module algebra for
    $\mathscr{U}=\overline{\mathscr{U}}_qs\ell(2)$ at an even root of
    unity}, Theor.\ Math. Phys. 159 (2009) 424--447
  [arXiv:\linebreak[0]0809.0144 [hep-th]].


\bibitem{[S-H]}A.M.\;Semikhatov, \textit{A Heisenberg double addition
    to the logarithmic Kazhdan--Lusztig duality},
  arXiv:\linebreak[0]0902.2215 [math.QA].


\bibitem{[NSz]} F.\;Nill and K.\;Szlach\'anyi, \textit{Quantum chains
    of Hopf algebras with quantum double cosymmetry}, Commun.\ Math.\
  Phys. 187 (1997) 159--200 [hep-th$/$\linebreak[0]9509100].


\bibitem{[CFM]}M.\;Cohen, D.\;Fischman, and S.\;Montgomery, \textit{On
    Yetter--Drin\-feld categories and $H$-commutativity}, Commun.\
  Algebra 27 (1999) 1321--1345.


\bibitem{[AGL]}A.\;Alekseev, D.\;Gluschenkov, and A.\;Lyakhovskaya,
  \textit{Regular representation of the quantum group $sl_q(2)$
    \textup{(}$q$ is a root of unity\textup{)}}, St.\ Petersburg
  Math. J.\ 6 (1994) 88.


\bibitem{[Su]}R.\;Suter, \textit{Modules over
    $\mathfrak{U}_q(\mathfrak{sl}_2)$}, Commun. Math. Phys. 163 (1994)
  359--393.


\bibitem{[X]}J.\;Xiao, \textit{Finite dimensional representations of
    $U_t(sl(2))$ at roots of unity}, Can. J. Math. 49 (1997) 772--787.


\bibitem{[Y]} D.N.\;Yetter, \textit{Quantum groups and representations
    of monoidal categories}, Math. Proc. Cambridge Philos.  Soc. 108
  (1990) 261--290.


\bibitem{[LR]}L.A.\;Lambe and D.E.\;Radford, \textit{Algebraic aspects
    of the quantum Yang--Baxter equation}, J. Alg.  154 (1992)
  228--288.


\bibitem{[RT]}D.E.\;Radford and J.\;Towber, \textit{Yetter--Drinfel'd
    categories associated to an arbitrary bialgebra}, J.\ Pure
  Appl. Algebra 87 (1993), 259--279.


\bibitem{[Mont]}S.\;Montgomery, \textsl{Hopf Algebras and Their
    Actions on Rings}, CBMS 82 (1993), American Mathematical Society,
  Providence, Rhode Island.


\bibitem{[Sch]}P.\;Schauenburg, \textit{Hopf modules and
    Yetter--Drinfel'd modules}, J. Algebra 169 (1994) 874--890.


\bibitem{[CGW]} M.\;Cohen, S.\;Gelaki, and S.\;Westreich, \textit{Hopf
    Algebras}, in: \textsl{Handbook of Algebra}, vol. 4.  Edited by
  M.\;Hazewinkel, Elsevier (2006) 173--239.


\bibitem{[CW]}M.\;Cohen and S.\;Westreich, \textit{From supersymmetry
    to quantum commutativity}, J.\ Algebra 168 (1994) 1--27.


\bibitem{[Mj]}S.\;Majid, \textit{$q$-Euclidean space and quantum group
    wick rotation by twisting}, J.\ Math.\ Phys. 35 (1994) 5025--5034
  [hep-th$/$\linebreak[0]9401112].


\bibitem{[AM]} D.\;Adamovi\'c and A.\;Milas, \textit{Lattice
    construction of logarithmic modules for certain vertex algebras},
  arXiv:\linebreak[0]0902.3417 [math.QA].


\bibitem{[NT]} K.\;Nagatomo and A.\;Tsuchiya, \textit{The triplet
    vertex operator algebra $W(p)$ and the restricted quantum group at
    root of unity}, arXiv:\linebreak[0]0902.4607 [math.QA].


\bibitem{[KS]}H.\;Kondo and Y.\;Saito, \textit{Indecomposable
    decomposition of tensor products of modules over the restricted
    quantum universal enveloping algebra associated to
    $\boldsymbol{\mathfrak{sl}_2}$}, arXiv:\linebreak[0]0901.4221
  [math.QA].


\bibitem{[Erd]} K.\;Erdmann, E.L.\;Green, N.\;Snashall, and
  R.\;Taillefer, \textit{Representation theory of the Drinfeld doubles
    of a family of Hopf algebras,} J.\ Pure and Applied Algebra 204
  (2006) 413--454 [math.RT$/$\linebreak[0]0410017].


\bibitem{[Sch-Galois]}P.\;Schauenburg, \textit{Galois objects over
    generalized Drinfeld doubles, with an application to
    $u_q(\mathfrak{sl}_2)$}, J.\ Algebra 217 (1999) 584--598


\bibitem{[FIK1]}A.T.\;Filippov, A.P.\;Isaev, and A.B.\;Kurdikov,
  \textit{Paragrassmann analysis and quantum groups},
  Mod. Phys. Lett. A7 (1992) 2129 [hep-th/9204089].

\end{thebibliography}
\end{document}